# New lower bounds for Schur and weak Schur numbers

Romain Ageron[a,1]   Paul Casteras[a,1]   Thibaut Pellerin[a,1]   Yann Portella[a,1]

Arpad Rimmel[b,2]   Joanna Tomasik[b,2,*]

[a]CentraleSupélec, Université Paris-Saclay, 91 190, Gif-sur-Yvette, France
[b]LISN, CentraleSupélec, Université Paris-Saclay, 91 405, Orsay, France

April 4, 2022

## Abstract

This article provides new lower bounds for both Schur and weak Schur numbers by exploiting a "template"-based approach. The concept of "template" is also generalized to weak Schur numbers. Finding new templates leads to explicit partitions improving lower bounds as well as the growth rate for Schur numbers, weak Schur numbers, and multicolor Ramsey numbers $R_n(3)$. The new lower bounds include $S(9) \geqslant 17\,803$, $S(10) \geqslant 60\,948$, $WS(6) \geqslant 646$, $WS(9) \geqslant 22\,536$ and $WS(10) \geqslant 71\,256$.

## 1 Introduction

We are interested in partitioning the set of integers $\{1, ..., p\}$ into $n$ subsets such that there is no subset containing three integers $x$, $y$ and $z$ verifying $x + y = z$. We say these subsets are *sum-free*. If we lessen the above constraint by only considering integers $x, y$ such that $x \neq y$, we say the subsets are *weakly sum-free*. The greatest $p$ for which there is a partition into $n$ sum-free subsets is called the $n^{\text{th}}$ Schur number and is denoted $S(n)$ [1]. Likewise for weakly sum-free partitions we define $WS(n)$ the $n^{\text{th}}$ weak Schur number [2]. Values of $S(n)$ and $WS(n)$ are known for small $n$ only.

### 1.1 State of the art

Up to recently, the most efficient generic construction for Schur numbers was given by Abbott and Hanson [3] in 1972 with a recursive construction. It gave the best lower bounds for all sufficiently large numbers. No equivalent was known for weak Schur numbers and, as a result, the best known partitions for large weak Schur numbers did not use the *weakly* sum-free property.

As for smaller weak Schur numbers, the best lower bounds were obtained by computations. Eliahou [4], Bouzy [5] and Rafilipojaona [6] improved the lower bounds with Monte-Carlo methods. This was the main approach during the past decade.

In 2020, Rowley introduced the notion of templates for Schur and Ramsey numbers [7] which generalizes Abbott and Hanson's construction and produces new lower bounds (and inequalities) for Schur numbers. He also provided two inequalities for weak Schur numbers [8]. Besides, these two inequalities do use the *weakly* sum-free property.

The hereunder tables recap the lower bounds for Schur and weak Schur numbers, respectively. Gray cells indicate the exact values.

---

[1]These authors contributed equally.
[2]These authors supervised the work of the main authors.
*Corresponding author (`joanna.tomasik@centralesupelec.fr`)

Table 1: Comparison of lower bounds for Schur numbers

| $n$ | 1 | 2 | 3 | 4 | 5 | 6 | 7 | 8 | 9 | 10 | 11 | 12 |
|---|---|---|---|---|---|---|---|---|---|---|---|---|
| State of the art | 1 | 4 | 13 | 44 | 160 [9] | 536 [10] | 1 696 [11] | 5 286 [7] | 17 694 [7] | 60 320 [7] | 201 696 [7] | 637 856 [7] |
| **Our results** | | | | | | | | | 17 803 | 60 948 | 203 828 | 644 628 |

Table 2: Comparison of lower bounds for weak Schur numbers

| $n$ | 1 | 2 | 3 | 4 | 5 | 6 | 7 | 8 | 9 | 10 | 11 | 12 |
|---|---|---|---|---|---|---|---|---|---|---|---|---|
| State of the art | 2 | 8 | 23 | 66 | 196 [12] | 642 [8] | 2 146 [8] | 6 976 [8] | 22 056 [8] | 70 778 [8] | 241 282 [8] | 806 786 [8] |
| **Our results** | | | | | | 646 | | | 22 536 | 71 256 | 243 794 | 815 314 |

## 1.2 Organization of this article

Our main contribution is a generalization of the concept of template to weak Schur numbers. Our templates provide new lower bounds (and inequalities) for weak Schur numbers. As a special case, our construction also includes a construction similar to Abbott and Hanson's [3], but for *weak* Schur numbers this time.

In Section 2, we explain Rowley's template-based construction for Schur numbers. Then, we give new templates, thus providing new lowers bounds and inequalities as well as showing that the growth rates for both Schur and Ramsey numbers $R_n(3)$ exceed 3.28.

In Section 3, we generalize the concept of template to weak Schur numbers and provide new lower bounds for weak Schur numbers. Then, we use a different approach and give a new lower bound for $WS(6)$.

Now, we introduce notations and definitions used throughout this article.

## 1.3 Definitions and notations

We start by defining sum-free and weakly sum-free subsets to introduce regular and weak Schur numbers. The set of positive natural numbers is denoted by $\mathbb{N}^* = \mathbb{N}\backslash\{0\}$.

**Definition 1.1.** *A subset $A$ of $\mathbb{N}$ is said to be sum-free if*

$$\forall (a,b) \in A^2, \ a+b \notin A.$$

**Definition 1.2.** *A subset $B$ of $\mathbb{N}$ is said to be weakly sum-free if*

$$\forall (a,b) \in B^2, \ a \neq b \Longrightarrow a+b \notin B.$$

Let us notice that a sum-free subset is also weakly sum-free, hence justifying the name of *weakly* sum-free subsets. Given $p$ and $n$ two integers, we are interested in partitioning the set of integers $\{1, 2, ..., p\}$, denoted by $[\![1,p]\!]$, into $n$ (weakly) sum-free subsets.

Schur proved [1] that given a number of subsets $n$, there is a value of $p$ such that $[\![1,q]\!]$ cannot be partioned into $n$ sum-free subsets for $q \geqslant p$. A similar property holds for weakly sum-free subsets [2]. These observations lead to the following definitions written for $n \in \mathbb{N}^*$:

**Definition 1.3.** *There is a largest integer denoted by $S(n)$ such that $[\![1, S(n)]\!]$ can be partitioned into $n$ sum-free subsets. $S(n)$ is called the $n^{th}$ Schur number.*

**Definition 1.4.** *There is a largest integer denoted by $WS(n)$ such that $[\![1, WS(n)]\!]$ can be partitioned into $n$ weakly sum-free subsets. $WS(n)$ is called the $n^{th}$ weak Schur number.*

Subsets of an $n$-partition are denoted $A_1, ..., A_n$. The smallest element of $A_i$ is denoted $m_i = \min(A_i)$. By ordering the subsets, we mean assuming that $m_1 < ... < m_n$. We will make clear in the text when we use a subset ordering.

A partition may be seen as a number coloring.

**Definition 1.5.** *The coloring associated to a partition $A_1, ..., A_n$ of $[\![1, p]\!]$ is the function $f$ such that $\forall x \in [\![1, p]\!], x \in A_{f(x)}$. Likewise, the partition associated to a coloring $f$ of $[\![1, p]\!]$ with $n$ colors is $\forall c \in [\![1, n]\!], A_c = f^{-1}(c)$.*

## 2 Templates for Schur numbers

We use Rowley's template-based constructions [7] for Schur numbers. To improve lower bounds for Schur and Ramsey numbers, he has introduced special sum-free partitions verifying some additional properties which can be extended using a method generalizing Abbott and Hanson's construction [3]. He called these partitions "templates", and we keep this nomination. Then, we propose new templates and use them to produce new lower bounds for Schur numbers.

### 2.1 Definition of $[\![1, S^+]\!]$

We begin by introducing *S-templates*, standing for *"Schur templates"*. The idea is to consider the first line of Figure 1 not as a combination of two blocs, like in [3] but as a whole, single construction. An S-template is then defined as a new object playing the role of the first line but with less, yet sufficient, constraints which allow for an expansion of the partition using a sum-free partition.

**Definition 2.1.** *Let $(p, n) \in (\mathbb{N}^*)^2$. An S-template with $n$ colors is defined as a partition of $[\![1, p]\!]$ into $n$ sum-free subsets $A_1, A_2, ..., A_n$ verifying for all subsets but $A_n$*

$$\forall (x, y) \in A_i^2, x + y > p \implies x + y - p \notin A_i. \tag{1}$$

*The integer $p$ is called the width of the template.*

Here $n$ can be seen as a "special" color in the sense that it is not constrained by Assertion (1). Number $p$ is necessarily colored with color $n$. However, note that $n$ is not necessarily the last color by order of appearance.

**Proposition 2.2.** *Let $n \in [\![2, +\infty[\![$. We define $S^+(n)$ as the greatest integer $w$ such that there is an S-template with width $w$ and $n$ colors. $S^+(n)$ is well defined and verifies*

$$2S(n-1) + 1 \leqslant S^+(n) \leqslant S(n).$$

*Proof.* The upper bound comes from the fact that an S-template with width $p$ and $n$ colors is also a sum-free $n$-partition of $[\![1, p]\!]$. The lower bound comes from Abbott and Hanson's construction [3]. □

### 2.2 Construction of Schur partitions using S-templates

We start by reminding the explicit construction of a sum-free partition with the use of an S-template and a sum-free partition. We rephrase in terms of Schur numbers this construction stated by Rowley in the context of Ramsey numbers [7].

**Theorem 2.3.** *Let $(p, k), (q, n) \in (\mathbb{N}^*)^2$. If there are an S-template with width $q$ and $n + 1$ colors and a sum-free $k$-partition of $[\![1, p]\!]$ then there is a partition of $[\![1, pq + m_{n+1} - 1]\!]$ into $n + k$ sum-free subsets. $m_{n+1}$ denotes the minimum element colored with the special color in the S-template.*

The idea lying beneath this theorem is similar to Abbott and Hanson's construction [3]. First, they put the integers from $[\![1, S(k)(2S(n) + 1) + S(n)]\!]$ in a table with width $2S(n) + 1$ and height $S(k) + 1$ and they numbered both the rows and the columns, starting from one upwards. Then, they extended a sum-free partition vertically by repeating it. Next, they used another sum-free partition to color the other half according to the line number, as in Figure 1. We give an example for $p = 4$, $q = 9$, $n = 2$ and $k = 2$ showing that $S(2+2) \geqslant S(2)(2S(2)+1) + S(2)$, both with Abbott and Hanson's construction (Figure 1) and with an S-template which is not included in Abbott and Hanson's construction (Figure 2). In both cases, the special color is grey.

Figure 1: Visualisation of an Abbott and Hanson construction

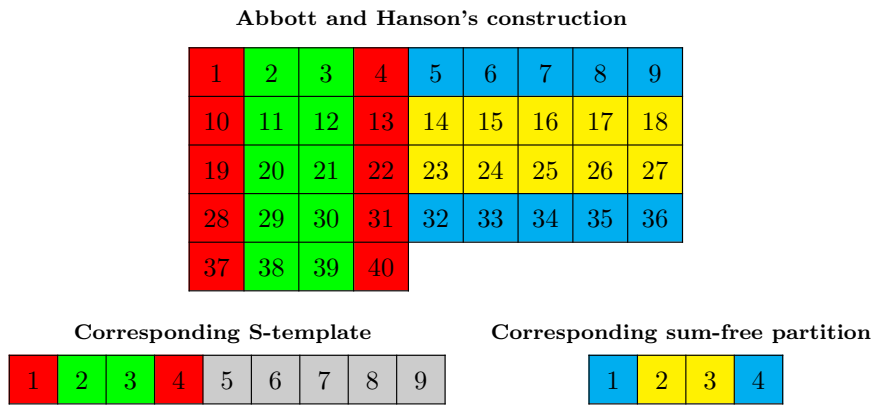

In the general construction with S-templates, the special color may not paint consecutive numbers any longer. However, the special color is still replaced by the colors of the sum-free partition according to the line number and the other colors are still vertically extended.

Figure 2: Visualisation of a template-based construction

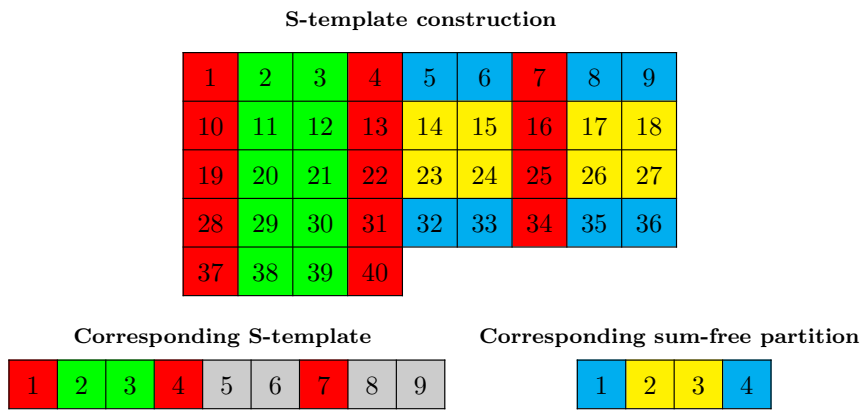

Now, we prove Theorem 2.3.

*Proof.* Let $f$ and $g$ be colorings associated to the S-template with width $q$ and the sum-free partition of $[\![1, p]\!]$, respectively: $f : [\![1, q]\!] \longrightarrow [\![1, n+1]\!]$ and $g : [\![1, p]\!] \longrightarrow [\![1, k]\!]$.

NB: In the following three predicates, the conditions $x + y \leqslant q$ and $x + y \leqslant p$ are omitted for readability.

The sum-free condition is expressed as:

$$\forall (x, y) \in [\![1, q]\!]^2, f(x) = f(y) \implies f(x+y) \neq f(x),$$

$$\forall (x, y) \in [\![1, p]\!]^2, g(x) = g(y) \implies g(x+y) \neq g(x).$$

The additional constraint for the S-template is:

$$\forall (x,y) \in [\![1,q]\!]^2, \begin{cases} f(x) = f(y) \leqslant n \\ x+y > q \end{cases} \implies f(x+y-q) \neq f(x).$$

For $x \in [\![1, pq+m_{n+1}-1]\!]$, we write $x = (\alpha-1)q+u$ for certain integers $\alpha \in \mathbb{N}^*$ and $u \in [\![1,q]\!]$. Obviously, this decomposition unique. Coefficients $\alpha$ and $u$ may be interpreted as line and column coordinates of $x$ in the template. We define a coloring $h$ of $[\![1, pq+m_{n+1}-1]\!]$ using $[\![1, n+k]\!]$ colors:

$$h: [\![1, pq+m_{n+1}-1]\!] \longrightarrow [\![1, n+k]\!],$$
$$x \longmapsto \begin{cases} f(u), & \text{if } f(u) \leqslant n, \\ n+g(\alpha), & \text{if } f(u) = n+1. \end{cases}$$

Function $h$ is well-defined since, by definition of $m_{n+1}$, $\forall x \in [\![pq+1, pq+m_{n+1}-1]\!]$, $f(u) \leqslant n$ and therefore $\forall x \in [\![1, pq+m_{n+1}-1]\!]$, $f(u) = n+1 \implies \alpha \in [\![1,p]\!]$.

Now, we prove that $h$ is a sum-free coloring. Let $x, y \in [\![1, pq+m_{n+1}-1]\!]$ such that $h(x) = h(y)$ and $x+y \leqslant pq+m_{n+1}-1$. We claim that $h(x+y) \neq h(x)$. We write $x = (\alpha-1)q+u$ and $y = (\beta-1)q+v$ where $\alpha, \beta \in \mathbb{Z}$ and $u, v \in [\![1,q]\!]$. Two cases are to be distinguished according to the value of $h(x)$.

**Case 1:** $h(x) \leqslant n$

Let us assume that $h(x+y) \leqslant n$, otherwise $h(x+y) \neq h(x)$ obviously holds. By definition of function $h$ and given that $h(u) = h(v)$, we conclude $f(u) = f(v)$. Two cases are to be distinguished according to the value of $x+y$.

- If $u+v > q$, we write $w = u+v-q \in [\![1,q]\!]$. Consequently $x+y = (\alpha+\beta-1)q+w$. By definition, $h(x+y) = f(w)$. Given that $f(u) = f(v) \leqslant n$, the additional constraint on $f$ implies $f(w) \neq f(u)$, that is $h(x+y) \neq h(x)$.

- If $u+v \leqslant q$, we write $w = u+v \in [\![1,q]\!]$. Consequently $x+y = (\alpha+\beta-2)q+w$. By definition, $h(x+y) = f(w)$. Given that $f(u) = f(v) \leqslant n$, the sum-free property of $f$ implies $f(w) \neq f(u)$, that is $h(x+y) \neq h(x)$.

**Case 2:** $h(x) \geqslant n+1$

Now we have $h(x) = n+g(\alpha) = n+g(\beta) = h(y)$, hence $g(\alpha) = g(\beta)$. As above, we distinguish between two cases according to the value of $x+y$.

- If $u+v > q$, write $w = u+v-q \in [\![1,q]\!]$. Then $x+y = (\alpha+\beta-1)q+w$. Assume that $h(x+y) \geqslant n+1$, otherwise $h(x+y) \neq h(x)$ obviously holds. By definition, $h(x+y) = n+g(\alpha+\beta)$. Given that $g(\alpha) = g(\beta)$, the sum-free property of $g$ implies $g(\alpha+\beta) \neq g(\alpha)$ that is $h(x+y) \neq h(x)$.

- If $u+v \leqslant q$, write $w = u+v \in [\![1,q]\!]$. Then $x+y = (\alpha+\beta-2)q+w$. The sum-free property of $f$ implies $f(w) \neq f(u)$. Therefore $f(w) \leqslant n$ and thus $h(x+y) \leqslant n$. In particular, given that $h(x) \geqslant n+1$, $h(x+y) \neq h(x)$.

$\square$

Setting $q = S^+(n+1)$ and $p = S(k)$ in Theorem 2.3 yields the following corollary.

**Corollary 2.4.** *Let $n, k \in \mathbb{N}^*$. Then*

$$S(n+k) \geqslant S^+(n+1)S(k) + m_{n+1} - 1.$$

The following proposition may improve the additive constant of Corollary 2.4.

**Proposition 2.5.** *Let $(q,n) \in (\mathbb{N}^*)^2$ and let $f$ be a coloring associated to an S-template with width $q$ and $n+1$ colors. Let $b \in \mathbb{N}$ and assume there is a coloring $g$ of $[\![1,b]\!]$ with $n+1$ colors such that:*

- $\forall (x,y) \in [\![1,q]\!]^2, \begin{cases} f(x) = f(y) \\ (x+y) \mod q \leqslant b \end{cases} \implies g((x+y) \mod q) \neq f(x),$

- $\forall (x,y) \in [\![1,q]\!] \times [\![1,b]\!], \begin{cases} f(x) = g(y) \\ x+y \leqslant b \end{cases} \implies g(x+y) \neq f(x).$

Then, for every $n \in \mathbb{N}^*$, by using on the last row the coloring $x \longmapsto g(x - qS(n))$, we have

$$S(n+k) \geqslant qS(k) + b.$$

This proposition corresponds to the fact that the hypotheses made on the coloring of the last row can be weakened because the constraints in the definition of an S-template prevent a number from:

- being the sum of two of numbers of its own color.
- creating a sum in its own color.

The later is not always relevant on the last row given that it is not completely filled. As a result, we can change the coloring on $[\![qS(k)+1, qS(k)+b]\!]$ (the last row) to extend the previous partitions.

There is a construction theorem for S-templates as well.

**Theorem 2.6.** *Let $(p,k),(q,n) \in (\mathbb{N}^*)^2$. If there are an S-template with width $q$ and $n+1$ colors, and an S-template with width $p$ and $k$ colors, then there also is an S-template with width $pq$ and $(n+k)$ colors.*

*Proof.* The idea is the same as in the Theorem 2.3. The only difference is the S-template property inherited from the second S-template.

$\square$

The inequality associated with this theorem is given by:

**Corollary 2.7.** *Let $n, k \in \mathbb{N}^*$. Then*

$$S^+(n+k) \geqslant S^+(n+1)S^+(k).$$

## 2.3 New lower bounds for Schur numbers

We now give the inequalities corresponding to the current best S-templates.

**Definition 2.8.** *A sum-free partition $A_1, ..., A_n$ of $[\![1,p]\!]$ is said to be symmetric if for all $x \in [\![1,p]\!]$, $x$ and $p+1-x$ belong to the same subset (except if $x = p+1-x$).*

*An S-template with $n$ colors is said to be symmetric if the sum-free $n$-partition derived from this template by applying the extension procedure of Theorem 2.3 with a sum-free partition with length one is symmetric.*

We produced S-templates using the lingeling SAT solver [13], hence providing lower bounds on $S^+$ and inequalities of the type $S(n+k) \geqslant aS(n) + b$. We sought templates providing us with the largest value of $(a,b)$ (in the lexicographic order). In order to reduce the search space when the number of colors exceeded five we only looked for symmetric S-templates, we assumed that the special color was the last color to appear and we constrained the $m_c$'s out of being too small. Details concerning the encoding as a SAT problem can be found in [9].

The following six inequalities are given by the current best S-templates with $n \leqslant 7$ colors.

$$S(n+1) \geqslant 3\,S(n) + 1 \tag{2}$$
$$S(n+2) \geqslant 9\,S(n) + 4 \tag{3}$$
$$S(n+3) \geqslant 33\,S(n) + 6 \tag{4}$$

Inequality (2) comes from Schur's original article [1]. Inequality (3) is due to Abott and Hanson [3] and inequality (4) to Rowley [7]. The three following inequalities are our result:

$$S(n+4) \geqslant 111\,S(n) + 43 \tag{5}$$
$$S(n+5) \geqslant 380\,S(n) + 148 \tag{6}$$
$$S(n+6) \geqslant 1\,160\,S(n) + 536 \tag{7}$$

The templates corresponding to inequalities (4), (5), and (6) are listed in A.

Inequalities (2), (3), and (4) cannot be further improved (with this definition of S-template). Inequality (5) cannot be further improved by only searching for symmetric S-templates whose special color is the last in the order of apparition (and with a multiplicative factor less than or equal to 118). Inequality (6) can most likely be further improved but the improvement probably will not be substantial. Finally, an S-template corresponding to inequality (7) was found by extending into a S-template the Schur 6-partition used in [11] (owing to the size of this template, its interest is limited and since it can easily be derived from above mentioned partition, it is not given in A). Although we could not find a better S-template with seven colors, inequality (7) can definitely be improved by a wide margin.

The previous inequalities give new lower bounds for $S(n)$ for $n \geqslant 9$. We compute the lower bounds for $n \in [\![8, 15]\!]$ using the four different inequalities. The best lower bounds are highlighted.

Table 3: New lower bounds for $S(n)$ with $n \in [\![8, 15]\!]$

| $n$ | 8 | 9 | 10 | 11 |
|---|---|---|---|---|
| $33\,S(n-3) + 6$ | **5 286** | 17 694 | 55 974 | 174 444 |
| $111\,S(n-4) + 43$ | 4 927 | **17 803** | 59 539 | 188 299 |
| $380\,S(n-5) + 148$ | 5 088 | 16 868 | **60 948** | **203 828** |

| $n$ | 12 | 13 | 14 | 15 |
|---|---|---|---|---|
| $33\,S(n-3) + 6$ | 587 505 | **2 011 290** | 6 726 330 | 21 272 730 |
| $111\,S(n-4) + 43$ | 586 789 | 1 976 176 | 6 765 271 | 22 624 951 |
| $380\,S(n-5) + 148$ | **644 628** | 2 008 828 | **6 765 288** | **23 160 388** |

Except for $S(8)$, $S(9)$, and $S(13)$, the best lower bounds are obtained thanks to the fifth inequality $S(n+5) \geqslant 380\,S(n) + 148$. The table does not go any further, but the same inequality allows one to improve the lower bounds for $n \geqslant 15$.

**Corollary 2.9.** *The growth rate for Schur numbers (and Ramsey numbers $R_n(3)$) satisfies $\gamma \geqslant \sqrt[5]{380} \approx 3.28$.*

*Proof.* It is a mere consequence of the inequality $S(n+5) \geqslant 380\,S(n) + 148$. As for Ramsey numbers, the following inequality holds $S(n) \leqslant R_n(3) - 2$ (see [3]) hence the result. □

## 2.4 Conclusion on S-templates

First, we formalized Rowley's template-based constructions [7] in the context of Schur numbers by introducing S-templates as well as a new sequence, $S^+$. We found new S-templates allowing us to obtain new lower bounds for schur numbers. One may notice that we mostly gave only lower bounds for $S^+$. It should be possible to find better S-templates by making different assumptions or using a different method (Monte-Carlo methods, for instance).

In the next section, we provide similar results for weak Schur numbers. We introduce WS-templates and a corresponding sequence, $WS^+$. Then, we derive similar relations and a construction method allowing us to find new lower bounds for weak Schur numbers.

# 3 Templates for weak Schur numbers

We generalize Rowley's constructions for weak Schur numbers [8]. We give an analogous for weak Schur numbers of Abbott and Hanson's construction for Schur numbers [3]. Then, as in the previous section, we introduce *WS-templates*, standing for *"weak Schur templates"*, as well as an associated sequence $WS^+(n)$.

We find templates and use them to provide new lower bounds for weak Schur numbers. Finally, we give a short explanation for the new lower bound $WS(6) \geqslant 646$ which was not directly obtained with a template, contrary to the other lower bounds given in this article.

## 3.1 Inequality for weak Schur numbers using Schur and weak Schur numbers

Up to now, no equivalent for weak Schur numbers of Abbott and Hanson's construction for Schur numbers [3] was known. Here we give a general lower bound for weak Schur numbers as a function of both regular and weak Schur numbers.

**Theorem 3.1.** *Let $(p,k), (q,n) \in (\mathbb{N}^*)^2$. If there are a weakly sum-free $n$-partition of $[\![1,q]\!]$ and a sum-free $k$-partition of $[\![1,p]\!]$ then there is a partition of $[\![1, p(q + \lceil \frac{q}{2} \rceil + 1) + q]\!]$ into $n+k$ weakly sum-free subsets.*

As Theorem 3.1 is actually a particular case of a more general theorem (Theorem 3.17) which will be formulated after the introduction of templates for weak Schur number, we only outline the demonstration. Its complete proof using templates for weak Schur numbers is in Subsection 3.3.

Let $(p,k), (q,n) \in (\mathbb{N}^*)^2$ such that there are a weakly sum-free $n$-partition of $[\![1,q]\!]$ and a sum-free $k$-partition of $[\![1,p]\!]$. Let $a \in \mathbb{N}$ with $a > q$ and let us try to build a coloring of $[\![1, ap+q]\!]$ into $n+k$ weakly sum-free subsets. Let $l = a - q - 1$, $r \in [\![1,q]\!]$ and $w = a - l - r - 1 = q - r$.

First, we construct the following table with $a$ columns and $p+1$ rows containing the integers from $[\![1, ap+q]\!]$. The columns are numbered from $-l$ to $+q$. This table is made of three different blocks.

- $\mathcal{T}$ (the "tail"): in the row indexed by 0: numbers from $[\![1,q]\!]$.

- $\mathcal{R}$ (the "rows"): the integers in columns from $-l$ to $+r$.

- $\mathcal{C}$ (the "columns"): the integers in the last $w$ columns.

Figure 3: Construction of the weakly sum-free coloring

| | | | | | | | 1 | 2 | ... | $r$ | $r+1$ | ... | $q-1$ | $q$ |
|---|---|---|---|---|---|---|---|---|---|---|---|---|---|---|
| | $a-l$ | $a-l+1$ | ... | $a-1$ | $a$ | $a+1$ | ... | $a+r-1$ | $a+r$ | $a+r+1$ | ... | $a+q-1$ | $a+q$ |
| | $2a-l$ | ... | ... | ... | $2a$ | ... | ... | ... | $2a+r$ | ... | ... | ... | $2a+q$ |
| $\mathcal{R}$ | ... | ... | ... | ... | ... | ... | ... | ... | ... | ... | ... | ... | ... |
| | ... | ... | ... | ... | ... | ... | ... | ... | ... | ... | ... | ... | ... |
| | $pa-l$ | ... | ... | ... | $pa$ | ... | ... | ... | $pa+r$ | ... | ... | ... | $pa+q$ |

$\mathcal{T}$ **block**
We color this block using the weakly sum-free coloring of $[\![1,q]\!]$ with colors $1, ..., n$.

$\mathcal{R}$ **block**
We use the colors $n+1, ..., n+k$ in this block. We color an integer $x$ according to its line number $\lambda(x)$. For any $x \in \mathcal{R}$, we color $x$ with $n+c$ where $c$ is the color of $\lambda(x)$ in the sum-free coloring of $[\![1,p]\!]$. Let $(x,y) \in \mathcal{R}^2$. The cases are twofold.

- $\lambda(x+y) = \lambda(x) + \lambda(y)$
  If $x + y \in \mathcal{R}$, we use the sum-free property of the coloring of $[\![1,p]\!]$. Otherwise, $x+y \in \mathcal{C}$ which is fine because we only use colors $1, ..., n$ in block $\mathcal{C}$.

- $\lambda(x+y) \neq \lambda(x) + \lambda(y)$
  In this case, we do not have information about the color of $\lambda(x+y)$. Thereby, we want to have $x+y \in \mathcal{C}$. A simple solution is to limit the horizontal displacement made by summing the two integers. That is, when the sum's row number is different from $\lambda(x) + \lambda(y)$, to upper bound $|x+y-(\lambda(x)+\lambda(y))a|$ so that $x+y$ stays in $\mathcal{C}$. There, the maximal displacement to the left (resp. to the right) is $2l$ (resp. $2r$). Not crossing entirely $\mathcal{C}$ by going to the left is then expressed as $-2l > -a+r$. Likewise, not going to far to the right is expressed as $2r < a - l$. It can then be written as $\max(l, r) \leqslant w$.

### $\mathcal{C}$ block

In this block, we use colors $1, ..., n$. We color an integer $x$ according to its projection on the first row $\pi(x)$ offseted by $a$, that is its column number. A simple solution is to color $x$ with the same color as $\pi(x) - a$ in the weakly sum-free coloring of $[\![1,q]\!]$. As long as $2q \leqslant a + r$ (not going two far to the right) and there is no $x \in \pi(\mathcal{C}) - a$ such that $2x \in \pi(\mathcal{C}) - a$, we do not have a sum in $\mathcal{C}$ when taking two numbers in the same colum. Finally, the weakly sum-free property of the coloring of $[\![1,q]\!]$ guarantees that if $x \in [\![1,r]\!]$ and $y \in \mathcal{C}$ are colored the same, then $x + y$ is colored differently.

Setting $w = l = \left\lceil \dfrac{q}{2} \right\rceil$ and $r = \left\lfloor \dfrac{q}{2} \right\rfloor$ satisfies these constraints: it corresponds to Theorem 3.1.

By setting $q = WS(n)$ and $p = S(k)$ in Theorem 3.1, one obtains the following corollary.

**Corollary 3.2.** $WS(n+k) \geqslant S(k)\left(WS(n) + \left\lceil \dfrac{WS(n)}{2} \right\rceil + 1\right) + WS(n)$ *for all* $(n,k) \in (\mathbb{N}^*)^2$.

This can be seen as an equivalent for weak Schur numbers of Abott and Hanson's construction for Schur numbers [3]. This formula also includes the results from [8] as a special case. For $n > 2$, the formula from Corollary 3.2 does not give new lower bounds.

**Remark 3.3.** *The inequality from Corollary 3.2 can be improved by adding 1 to the lower bound if $WS(n)$ is odd (more generally if $q$ is odd in Theorem 3.1). However, it lengthens the proof and it is never useful in practice.*

Now, we introduce WS-templates and the sequence $WS^+$ in order to generalize the above construction.

## 3.2 Definition of $[\![1, WS^+]\!]$

In this subsection, we introduce WS-templates and prove calculative results for the general construction theorem on templates for weak Schur numbers.

**Definition 3.4.** *Let $(a,b) \in (\mathbb{N}^*)^2$ such that $a > b$. We define:*

$$\pi_{a,b} : x \longmapsto (\mathrm{Id} + a\mathbb{1}_{[\![0,b]\!]})(x \mod a).$$

If there is no confusion on the $a$ and $b$ to use, $\pi_{a,b}$ is denoted by $\pi$. Notice that for all $x \in \mathbb{Z}$, $\pi(x) = x \mod a$ and for all $x \in [\![b+1, a+b]\!], b+1 \leqslant \pi(x) \leqslant a+b$.

Function $\pi$ is the projection on the first line mentioned in the intuitive explanation. The following four propositions are calculative properties on $\pi$ reflecting the behaviour of an element's column number in Figure 3 and that we will use later when we introduce WS-templates.

**Proposition 3.5.**
$$\forall x \in [\![b+1, a+b]\!], \pi(x) = x.$$

*Proof.* Let $x \in [\![b+1, a+b]\!]$. If $x < a$ then $x \mod a = x \notin [\![1,b]\!]$. Hence $\pi(x) = x$. Otherwise, $x \mod a = x - a \in [\![1,b]\!]$. Hence $\pi(x) = x - a + a = x$. □

**Proposition 3.6.** *Let $x \in [\![1,b]\!]$ and $y \in \mathbb{N}^*$. Then*

$$\pi(x + \pi(y)) = \pi(x + y).$$

*Proof.* It is a direct consequence of $\pi(x) = x \mod a$.

□

**Proposition 3.7.** *Let $x \in [\![1, b]\!]$ and $y \in \mathbb{N}^*$ such that $x + \pi(y) \leqslant a + b$. Then*

$$\pi(x + y) = x + \pi(y).$$

*Proof.* $\pi(y) \geqslant b + 1$ and thus $b + 1 \leqslant x + \pi(y) \leqslant a + b$.

$$\begin{aligned} \pi(x + y) &= \pi(x + \pi(y)) \quad \text{by Proposition 3.6} \\ &= x + \pi(y) \quad \text{by Proposition 3.5} \end{aligned}$$

□

**Proposition 3.8.** *Let $(x, y) \in (\mathbb{N}^*)^2$. Then*

$$\pi(\pi(x) + \pi(y)) = \pi(x + y).$$

*Proof.* It is a direct consequence of $\pi(x) = x \mod a$.

□

After defining the function related to the column number of each element in Figure 3, we introduce the function related to its line number.

**Definition 3.9.** *Let $(a, b) \in (\mathbb{N}^*)^2$ such that $a > b$. Define*

$$\lambda_{a,b} : x \longmapsto 1 + \left\lfloor \frac{x - b - 1}{a} \right\rfloor.$$

If there is no confusion on the $a$ and $b$ to use, $\lambda_{a,b}$ is denoted by $\lambda$.

Function $\lambda$ maps an element $x$ to its line number as mentioned in the intuitive explanation. As we have just done with $\pi$, we prove three more calculative results on both $\pi$ and $\lambda$ that are used in Subsection 3.3.

**Proposition 3.10.** *Let $x \in \mathbb{N}^*$. Then*

$$x = a\lambda(x) + \pi(x) - a.$$

*Proof.* Let $(a, b) \in (\mathbb{N}^*)^2$ such that $a > b$ and let $x \in \mathbb{N}^*$.

We have $a\lambda(x) + \pi(x) - a = a \left\lfloor \frac{x - b - 1}{a} \right\rfloor + (x \mod a) + a\mathbb{1}_{[\![0,b]\!]}(x \mod a)$.

- If $x \mod a > b$ then $a\lambda(x) + \pi(x) - a = a \left\lfloor \frac{x}{a} \right\rfloor + x \mod a = x$.

- If $x \mod a \leqslant b$ then $a\lambda(x) + \pi(x) - a = a \left( \left\lfloor \frac{x}{a} \right\rfloor - 1 \right) + x \mod a + a = x$.

□

**Proposition 3.11.** *Let $x, y \in \mathbb{Z}$ such that $\lambda(x + y) = \lambda(y)$. Then*

$$\pi(x + y) = x + \pi(y).$$

*Proof.* By applying Proposition 3.10 twice, we get $a\lambda(x+y) + \pi(x+y) - a = x + y = x + a\lambda(y) + \pi(y) - a$. The result is then obtained by simplifying the equality.

□

**Proposition 3.12.** *Let $x, y \in \mathbb{Z}$ such that $\pi(x) + \pi(y) \in [\![a + b + 1, 2a + b]\!]$. Then*

$$\lambda(x + y) = \lambda(x) + \lambda(y).$$

*Proof.* By Proposition 3.10, $x + y = a(\lambda(x) + \lambda(y)) + \pi(x) + \pi(y) - 2a$. Then

$$
\begin{aligned}
\lambda(x+y) &= \left\lfloor \frac{x+y-b-1}{a} \right\rfloor + 1 \\
&= \left\lfloor \frac{a(\lambda(x) + \lambda(y)) + \pi(x) + \pi(y) - 2a - b - 1}{a} \right\rfloor + 1 \\
&= \lambda(x) + \lambda(y) - 1 + \left\lfloor \frac{\pi(x) + \pi(y) - b - 1}{a} \right\rfloor \\
&= \lambda(x) + \lambda(y) - 1 + 1 \quad \text{since } \pi(x) + \pi(y) \in [\![a+b+1, 2a+b]\!] \\
&= \lambda(x) + \lambda(y).
\end{aligned}
$$

□

**Definition 3.13.** *Let $(a, n, b) \in (\mathbb{N}^*)^3$ with $a > b$. Let $(A_1, ..., A_n)$ a partition of $[\![1, a + b]\!]$. This partition is said to be a b-weakly-sum-free template (b-WS-template) with width $a$ and $n$ colors when:*

- $\forall i \in [\![1, n]\!]$, $A_i$ *is weakly-sum-free,*
- $\forall i \in [\![1, n]\!]$, $A_i \backslash [\![1, b]\!]$ *is sum-free,*
- *For $A_n$ (the special subset):*

$$\forall (x, y) \in A_n^2, \; x + y > b + 2a \implies x + y - 2a \notin A_n,$$

- *For the others subsets*

$$\forall i \in [\![1, n-1]\!], \; \forall (x, y) \in A_i^2, \; x + y > a + b \implies \pi(x + y) \notin A_i.$$

Number $a$ is necessarily colored with color $n$. Note that the special color $n$ is not necessarily the last color by order of appearance. Now, we introduce the number $WS^+(n)$ that plays the same role as $S^+(n)$ for S-templates.

**Definition 3.14.** *Let $(n, b) \in (\mathbb{N}^*)^2$. If there is $a$ such that there is a b-WS-template with width $a$ and $n$ colors, we define:*

$$WS_b^+(n) = \max\{a \in \mathbb{N}^* / \text{there is a b-WS-template with width } a \text{ and } n \text{ colors}\}.$$

*If there is no such $a$, we set $WS_b^+(n) = 0$.*

**Definition 3.15.** *Let $n \in \mathbb{N}^*$. We define:*

$$WS^+(n) = \max_{b \in \mathbb{N}^*} WS_b^+(n).$$

The following proposition shows how $WS^+(n)$ compares to weak Schur numbers.

**Proposition 3.16.** *Let $n \in [\![2, +\infty]\!]$. Then:*

$$\frac{3}{2} WS(n-1) + 1 \leqslant WS^+(n) \leqslant WS(n).$$

*Proof.* The lower bound comes from Corollary 3.2. The upper bound comes from the fact that a WS-template with width $a$ and $n$ colors is in particular a partition of $[\![1, a]\!]$ into $n$ weakly sum-free subsets. □

Having established these properties, we prove our main result.

## 3.3 Construction of weak Schur partitions using WS-templates

Theorem 2.3 describes a manner to extend an S-template into a larger Schur partitions using other Schur partitions. Likewise, a WS-template can be extended into larger weak Schur partitions by using Schur partitions. This is the object of the following theorem.

**Theorem 3.17.** *Let $(a, n, b) \in (\mathbb{N}^*)^3$ with $a > b$ and $(p, k) \in (\mathbb{N}^*)^2$. If there are a sum-free k-partition of $[\![1, p]\!]$ and a b-WS-template $(A_1, ..., A_{n+1})$ with width $a$ and $n + 1$ colors, then there is a partition of $[\![1, pa + b]\!]$ into $k + n$ weakly sum-free subsets.*

*Proof.* Let $(a, n, b) \in (\mathbb{N}^*)^3$ and $(p, k) \in (\mathbb{N}^*)^2$. Let $f$ and $g$ be colorings associated to the $b$-WS-template and to the sum-free partition of $[\![1, p]\!]$, respectively: $f : [\![1, a + b]\!] \longrightarrow [\![1, n + 1]\!]$ and $g : [\![1, p]\!] \longrightarrow [\![1, k]\!]$. Moreover, we assume that the sum-free coloring of $[\![1, p]\!]$ is ordered.

NB: To keep the notation short, the conditions $x + y \leqslant p$ and $x + y \leqslant a + b$ are omitted in the following five predicates.

The sum-free conditions are expressed as:
$$\forall (x, y) \in [\![b + 1, a + b]\!]^2, f(x) = f(y) \implies f(x + y) \neq f(x),$$
$$\forall (x, y) \in [\![1, p]\!]^2, g(x) = g(y) \implies g(x + y) \neq g(x).$$

The weakly sum-free condition is expressed as:
$$\forall (x, y) \in [\![1, a + b]\!]^2, \begin{cases} f(x) = f(y) \\ x \neq y \end{cases} \implies f(x + y) \neq f(x),$$

The additional constraints for the WS-template are expressed as:
$$\forall (x, y) \in [\![1, a + b]\!]^2, \begin{cases} f(x) = f(y) \leqslant n \\ x + y > a + b \end{cases} \implies f(\pi(x + y)) \neq f(x),$$
$$\forall (x, y) \in [\![1, a + b]\!]^2, \begin{cases} f(x) = f(y) = n + 1 \\ x + y > 2a + b \end{cases} \implies f(x + y - 2a) \neq f(x).$$

We split $[\![1, pa + b]\!]$ into three subsets.

NB: To keep the notation short, the restriction to $[\![b + 1, pa + b]\!]$ of $\pi$ defined in Subsection 3.2 is denoted by $\pi$ in the hereunder equations:

- $\mathcal{T} = [\![1, b]\!]$,
- $\mathcal{C} = \pi^{-1}(f^{-1}([\![1, n]\!]))$,
- $\mathcal{R} = \pi^{-1}(f^{-1}(\{n + 1\}))$.

We define a coloring $h$ of $[\![1, pa + b]\!]$ using $[\![1, n + k]\!]$ colors:
$$h : \begin{array}{rcl} [\![1, pa + b]\!] & \longrightarrow & [\![1, n + k]\!], \\ x & \longmapsto & \begin{cases} f(x), & \text{if } x \in \mathcal{T}, \\ f(\pi(x)), & \text{if } x \in \mathcal{C}, \\ n + g(\lambda(x)), & \text{if } x \in \mathcal{R}. \end{cases} \end{array}$$

Function $h$ is well defined since sets $\mathcal{T}, \mathcal{C}, \mathcal{R}$ partition $[\![1, pa + b]\!]$. Now, we prove that $h$ is a weakly sum-free coloring. Let $x, y \in [\![1, pa + b]\!]$ be such that $x \neq y$, $h(x) = h(y)$, and $x + y \leqslant pa + b$. We claim that $h(x + y) \neq h(x)$. Nine cases are to be distinguished according to which of the sets $\mathcal{T}, \mathcal{C}, \mathcal{R}$ $x$ and $y$ belong to. It is sufficient to check only six cases out of nine since $x$ and $y$ play symmetric roles.

**Case 1:** $(x, y) \in \mathcal{T}^2$

If $x + y \leqslant b$ then $h(x + y) = f(x + y)$. Otherwise, $b < x + y < a + b$ since $b < a$ and therefore $\pi(x + y) = x + y$ (Proposition 3.5). Hence in both cases $h(x + y) = f(x + y)$. Given that $f$ is a weakly sum-free coloring, $f(x + y) \neq f(x)$ since $f(x) = h(x) = h(y) = f(y)$ and $x \neq y$. That is $h(x + y) \neq h(x)$.

**Case 2:** $(x, y) \in \mathcal{T} \times \mathcal{C}$

Given that $h(x) = h(y)$ and by definition of $h$, $f(x) = f(\pi(y))$. Besides, $f(\pi(y)) \leqslant n$ since $y \in \mathcal{C}$. Two cases are to be distinguished according to the value of $x + \pi(y)$.

- If $x+\pi(y) \leqslant a+b$ then $f(x+\pi(y)) = f(\pi(x+y))$ (Proposition 3.7). Given that $f$ is a weakly sum-free coloring, $f(x+\pi(y)) \neq f(x)$ since $f(x) = f(\pi(y))$ and $x \neq \pi(y)$ since $x \leqslant b < \pi(y)$.

- If $x + \pi(y) > a + b$ then given that $f$ is a WS-template and since $f(x) = f(\pi(y)) \leqslant n$, $f(\pi(x+\pi(y))) \neq f(x)$. Furthermore $f(\pi(x+\pi(y))) = f(\pi(x+y))$ (Proposition 3.6), such that $f(\pi(x+y)) \neq f(x)$.

In both cases $f(\pi(x+y)) \neq f(x)$. If $f(\pi(x+y)) \leqslant n$ then $h(x+y) = f(\pi(x+y))$. Therefore $h(x+y) \neq h(x)$ since $f(x) = h(x)$. Otherwise, $f(\pi(x+y)) = n+1$ and thus $h(x+y) > n$. In particular, $h(x+y) \neq h(x)$ since $h(x) = h(y) \leqslant n$.

**Case 3:** $(x,y) \in \mathcal{T} \times \mathcal{R}$

Necessarily $h(x) = h(y) = n+1$. Two cases are to be distinguished according to the value of $\lambda(x+y)$.

- If $\lambda(y) = \lambda(x+y)$ then $\pi(x+y) = x + \pi(y)$ (Proposition 3.11). By definition of $h$, $f(x) = f(\pi(y))$. Given that $f$ is a weakly sum-free coloring, $f(x+\pi(y)) \neq f(x)$ since $f(x) = f(\pi(y))$ and $x \neq \pi(y)$ since $x \leqslant b < \pi(y)$. Hence $h(x+y) \neq h(x)$.

- If $\lambda(y) \neq \lambda(x+y)$ then $\lambda(x+y) = \lambda(y) + 1$ since $x \leqslant b < a$. Besides, $n+1 = h(y) = n + g(\lambda(y))$. Hence $g(\lambda(y)) = 1$. Moreover $g(1) = 1$ since $g$ is an ordered coloring. Therefore, given that $g$ is sum-free, $g(\lambda(y)+1) \neq 1$. If $\pi(x+y) \in A_{n+1}$ then $h(x+y) = n + g(\lambda(x+y)) \neq n+1$. Otherwise, $h(x+y) \leqslant n$.

In both cases $h(x+y) \neq h(x)$.

**Case 4:** $(x,y) \in \mathcal{C}^2$

By definition of $h$ and since $h(x) = h(y)$, $f(\pi(x)) = f(\pi(y))$. Two cases are to be distinguished according to the value of $\pi(x) + \pi(y)$.

- If $\pi(x) + \pi(y) \leqslant a+b$ then $\pi(x) + \pi(y) = \pi(x+y)$. Hence $f(\pi(x+y)) \neq f(\pi(x))$ since $f$ is sum-free for $x > b$

- If $\pi(x) + \pi(y) > a + b$ then given that $f$ is a WS-template, $f(\pi(\pi(x)+\pi(y))) \neq f(\pi(x))$ since $f(\pi(x)) = f(\pi(y))$. Besides, $f(\pi(\pi(x)+\pi(y))) = f(\pi(x+y))$ (Proposition 3.8). Hence $f(\pi(x+y)) \neq f(\pi(x))$.

In both cases $f(\pi(x+y)) \neq f(x)$. If $f(\pi(x+y)) \leqslant n$ then $h(x+y) = f(\pi(x+y))$. Therefore $h(x+y) \neq h(x)$ since $f(x) = h(x)$. Otherwise, $f(\pi(x+y)) = n+1$ and thus $h(x+y) > n$. In particular, $h(x+y) \neq h(x)$ since $h(x) = h(y) \leqslant n$.

**Case 5:** $(x,y) \in \mathcal{C} \times \mathcal{R}$

By definition of function $h$, $h(x) \neq h(y)$.

**Case 6:** $(x,y) \in \mathcal{R}^2$

In particular $f(\pi(x)) = f(\pi(y)) = n+1$. Three cases are to be distinguished according to the value of $\pi(x) + \pi(y)$.

- If $\pi(x) + \pi(y) \in [\![a+b+1, 2a+b]\!]$ then $\lambda(x+y) = \lambda(x) + \lambda(y)$ (Proposition 3.12). By definition of $h$ and since $h(x) = h(y)$, $g(\lambda(x)) = g(\lambda(y))$. Hence, $h(\lambda(x+y)) \neq h(\lambda(x))$ since $h$ is a sum-free coloring. If $f(x+y) \geqslant n+1$ then $h(x+y) = n + g(\lambda(x+y))$. And $h(x) = n + g(\lambda(x))$. Therefore, $h(x+y) \neq h(x)$. Otherwise $h(x+y) \leqslant n < h(x)$. In particular $h(x+y) \neq h(x)$.

- If $\pi(x) + \pi(y) > 2a + b$ then $f(\pi(\pi(x)+\pi(y))) \neq f(\pi(x)) = n+1$ since $f$ is a $b$-WS template and $f(\pi(x)) = f(\pi(y))$. Given that $\pi(\pi(x)+\pi(y)) = \pi(x+y)$ (Proposition 3.8), $f(\pi(x+y)) \neq n+1$.

- If $\pi(x) + \pi(y) \leqslant b + a$ then, given that $\pi(x) + \pi(y) \geqslant b$ and $f_{|[\![b,a+b]\!]}$ is sum-free, $f(\pi(x)+\pi(y)) \neq f(\pi(x)) = n+1$. That is $f(\pi(x+y)) \neq n+1$ (Proposition 3.5).

In both of the last two cases, $f(\pi(x+y)) \neq n+1$ that is $x+y \in \mathcal{C}$. Therefore $h(x+y) < n \leqslant h(x)$. In particular, $h(x+y) \neq h(x)$.

□

We go back to Theorem 3.1 which can be seen as a particular case of Theorem 3.17 in the same way Abott and Hanson's construction [3] can be seen as a particular case of Theorem 2.3. Its deferred proof is given below.

PROOF OF THEOREM 3.1. Let $(q, n) \in (\mathbb{N}^*)^2$ such that there is a partition of $[\![1, q]\!]$ into n weakly sum-free subsets. Let $f : [\![1, q]\!] \to [\![1, n]\!]$ be a weakly sum-free colouring. Let $b = q$ and $a = q + \left\lceil \frac{q}{2} \right\rceil + 1$. A new colouring $g$ is defined as follows:

$$g : \quad [\![1, a+b]\!] \longrightarrow [\![1, n+1]\!],$$
$$x \longmapsto \begin{cases} f(x), & \text{if } x \in [\![1, b]\!], \\ n + 1, & \text{if } x \in [\![b+1, 2b+1]\!], \\ f(x - a), & \text{if } x \in [\![2b+2, a+b]\!]. \end{cases}$$

We claim that $g$ is a $b$-WSF-template with width $a$ and $n + 1$ colours.

- Function $g_{|[\![b+1,a+b]\!]}$ is a sum-free colouring. Indeed, let $(x, y) \in [\![b+1, a+b]\!]^2$ such that $g(x) = g(y)$. If $g(x) = n + 1$ then $z = x + y > 2b + 1$ and therefore, either $z > a + b$ or $f(z) \neq n + 1$. Otherwise, $x + y > a + b$.

- Function $g$ is a weakly sum-free colouring. Indeed, let $(x, y) \in [\![1, a+b]\!]^2$ such that $z = x + y \leqslant a + b$ and $g(x) = g(y)$. Given that $x$ and $y$ have symmetric roles, we can assume that $x \leqslant y$. If $x > b$ then $g(z) \neq g(x)$ as seen above. If $y \leqslant b$ then $f(x) = g(x) = g(y) = f(y) \leqslant n$ and either $z \leqslant a + b$ and $f(z) \neq f(x)$ since $f$ is a weakly sum-free colouring or $a + b + 1 \leqslant z \leqslant 2b$ and $g(z) = n + 1$; therefore $g(z) \neq g(x)$. If $x \leqslant b$ and $y > b$ then $g(x) = f(x)$, $g(y) = f(y - a)$ and $g(z) = f(z - a)$. We have $x \neq y - a$ (otherwise, we would have $a + b \geqslant z = 2y - a \geqslant 4b + 4 - a > a + b$) and thus $f(x + y - a) \neq f(x)$, that is $g(z) \neq g(x)$.

- Colour $n + 1$ verifies the additional constraints for the special colour. Indeed, $b + 2a \geqslant 4b + 2$. Hence, $\forall (x, y) \in g^{-1}(\{n + 1\}), x + y \leqslant b + 2a$.

- Colours $1, ..., n$ verify the additional constraints for the regular colours. Indeed, let $(x, y) \in g^{-1}([\![1, n]\!])$ such that $x + y > a + b$. Given that $x$ and $y$ have symmetric roles, we can assume that $x \leqslant y$. Necessarily $y \geqslant 2b + 2$. If $x \leqslant b$ then $z = x + y \in [\![a + b + 1, a + 2b]\!]$ and therefore $\pi(z) \in [\![b + 1, 2b]\!]$. Otherwise, $z = x + y \in [\![4b + 4, 2a + 2b]\!]$ and therefore $\pi(z) \in [\![b + 1, 2b]\!]$. Hence, in both cases, $g(\pi(z)) \neq g(x)$.

The result is then obtained by applying Theorem 3.17.

$\square$

By setting $p = S(k)$ and $a = WS^+(n + 1)$ in Theorem 3.17, one obtains:

**Corollary 3.18.** *Let $n, k \in \mathbb{N}^*$ and set $b_{max} = \max\{b \in \mathbb{N}^* / WS_b^+(n + 1) = WS^+(n + 1)\}$. Then:*

$$WS(n + k) \geqslant S(k) WS^+(n + 1) + b_{max}.$$

**Remark 3.19.** *In the S-template construction for Schur numbers, the additive constant comes from the fact that the special color does not necessarily appear at the very begining of the repeating pattern. Likewise, $b_{max}$ can actually be replaced by*

$$\max_{b \in \mathbb{N}^*} \left\{ \min(A_{n+1} \backslash [\![1, b]\!]) - 1 \mid WS_b^+(n + 1) = WS^+(n + 1) \right\}.$$

As in Corollary 2.4, the additive constant Theorem 3.17 can be improved by weakening the hypotheses made on the last row. The principle behind it is the same as in Proposition 2.5.

**Proposition 3.20.** *Let $(b, k, a) \in (\mathbb{N}^*)^3$ and let $f$ be a coloring associated to a $b$-WS-template with width $p$ and $k$ colors. Let $c \in \mathbb{N}$ and assume there is a coloring $g$ of $[\![b + 1, b + c]\!]$ with $k$ colors such that for all $c \in [\![1, k]\!]$,*

- $\forall (x, y) \in [\![1, a+b]\!] \times [\![b+1, a+b]\!], \begin{cases} f(x) = f(y) \\ \pi(x + y) \leqslant b + c \end{cases} \implies g(\pi(x + y)) \neq f(x),$

- $\forall (x,y) \in [\![1, a+b]\!] \times [\![b+1, b+c]\!], \begin{cases} f(x) = g(y) \\ \pi(x+y) \leqslant b+c \end{cases} \implies g(\pi(x+y)) \neq f(x).$

Then, for every $n \in \mathbb{N}^*$, by using on the last row the coloring $x \longmapsto g(x - pS(n))$, we have

$$WS(n+k) \geqslant WS^+(k+1)S(n) + b + c.$$

The WS-templates can actually be fine-tuned further. However, it only gives minor improvements (most likely only an additive constant) while increasing the size of the search space.

These modifications rely on one main observation: the first row (excluding the "tail") has constraints that other rows do not have because of the tail, especially if the special color appears in the tail as well. Indeed, the line number of the sum of two numbers in the tail is never greater than one. Therefore, having a coloring dedicated to the first row and another coloring dedicated to the repeating pattern would weaken the constraints on the repeating pattern.

These additional constraints on the first row imply additional constraints on the rows of the its color (compared to other colors). However, one may notice that there are never two consecutive rows of this color. Moreover these additional constraints mostly concern two consecutive rows. Having a special coloring for rows of this color as well as another color dedicated to the first row avoids enforcing these constraints specific to one color on the other colors and thus on the main repeating pattern.

Finally, given that the color of the first row is different from the special color of the main repeating pattern, coloring $g$ in Proposition 3.20 can exploit this difference.

There is a construction theorem for WS-templates as well.

**Theorem 3.21.** *Let $(k, p) \in (\mathbb{N}^*)^2$ and $(a, n, b) \in (\mathbb{N}^*)^3$. If there are an S-template with width $p$ and $k+1$ colors and a b-WS-template with width $a$ and $n$ colors, then there is a pb-WS-template with width $pq$ and $(n+k)$ colors.*

*Proof.* The idea is the same as in theorem 3.17. The properties specific to WS-templates are inherited from both the S-template and the WS-template. $\square$

Theorem 3.21 yields the following corollary.

**Corollary 3.22.** *Let $n, k \in \mathbb{N}^*$. Then*

$$WS^+(n+k) \geqslant S^+(k+1) WS^+(n).$$

## 3.4 New lower bounds for Weak Schur numbers

We produced WS-templates using the lingeling SAT solver [13], hence providing lower bounds on $WS^+$ and inequalities of the type $WS(n+k) \geqslant aS(n) + b$. We sought templates providing the greatest value of $(a, b)$ (in the lexicographic order). Details concerning the encoding as a SAT problem can be found in [9].

Here are the inequalities given by the current best WS-templates. The template corresponding to the third inequality can be found in B.

$$WS(n+1) \geqslant 4\,S(n) + 2 \tag{8}$$
$$WS(n+2) \geqslant 13\,S(n) + 8 \tag{9}$$

Inequalities (8) and (9) were found by Rowley, they are detailed in [8].

$$WS(n+3) \geqslant 42\,S(n) + 24 \tag{10}$$
$$WS(n+4) \geqslant 132\,S(n) + 26 \tag{11}$$

Inequality (10), found with a SAT solver [13], cannot be improved (with this definition of WS-template). It uses the first sophistication explained in Subsection 3.3 in order to add the last number in the first color. As for inequality (11), it was obtained by combining an S-template with width 33 with a WS-template with width 4. The best template we could find with a computer search gives the inequality $WS(n + 4) \geqslant 127 S(n) + 68$. It was also found with the SAT solver. In order to reduce the search space, we only looked for WS-templates with five colors which start with a near-optimal $WS(4)$ partition and we assumed that the special color was the last by order of appearance. We expect that better WS-templates with $n \geqslant 5$ colors can be found but one would have not to use the above assumptions.

Like in Subsection 2.3, we compute the lower bounds given by Inequalities (8), (9), and (10) for $n \in [\![8, 15]\!]$. The best lower bound for each value of $n$ is highlighted.

Table 4: New lower bounds for $WS(n)$ with $n \in [\![8, 15]\!]$

| $n$ | 8 | 9 | 10 | 11 |
|---|---|---|---|---|
| $4\,S(n-1) + 2$ | 6 786 | 21 146 | 71 214 | **243 794** |
| $13\,S(n-2) + 8$ | **6 976** | 22 056 | 68 726 | 231 447 |
| $42\,S(n-3) + 24$ | 6 744 | **22 536** | **71 256** | 222 036 |

| $n$ | 12 | 13 | 14 | 15 |
|---|---|---|---|---|
| $4\,S(n-1) + 2$ | **815 314** | 2 578 514 | 8 045 162 | 27 061 154 |
| $13\,S(n-2) + 8$ | 792 332 | **2 649 772** | 8 380 172 | 26 146 778 |
| $42\,S(n-3) + 24$ | 747 750 | 2 559 840 | **8 560 800** | **27 074 400** |

## 3.5 Another approach: $WS(6) \geqslant 646$

The weak Schur partitions obtained with WS-templates have an extremely regular structure. Therefore, one may expect to find larger partitions by allowing more freedom in the structure of the partitions while still preserving a structure somewhat close to a template-based partition thus reducing the search space to a manageable size.

When applying the inequality $WS(n + 1) \geqslant 4\,S(n) + 2$, one may realize that it is possible to build a weakly sum-free partition of length $4\,S(n) + 3$ for small values of $n$ ($n \leqslant 4$) by using the construction of Corollary 3.18 for the integers $1, 2, 4i$ and $4i + 1$ for $i \in [\![1, S(n)]\!]$ but not constraining the other integers. We did the same for the Schur partitions corresponding to Schur number five. More precisely, we imposed these constraints only for $i \leqslant 50$ in order to have more degrees of freedom. The number 50 was chosen arbitrarly so that all of the Schur number 5-partitions could be tested in a few hours.

However, trying out all of the 2 447 113 088 Schur number 5-partitions [9] one by one would not result in a reasonable computation time: it is necessary to test the construction on several partitions at once. We encoded the problem as a satisfiability problem and used the 1616 backdoors that were used in [9] in order to encode a group of partitions in a compact and efficient way. Among all of the 1616 backdoors, only the 911[th] backdoor gave a weakly sum-free partition of length 643. This backdoor gave a weakly sum-free partition of length 646 and cannot give a weakly sum-free partition of length 647. A weakly sum-free partition of length 646 can be found in C.

## 3.6 Conclusion on WS-templates

We started by giving a new construction which can be seen as an equivalent for weakly sum-free partitions of Abott and Hanson's construction for sum-free partitions. We then generalized this construction by introducing WS-templates. This allows us to find new lower bounds and new inequalities for weak Schur numbers. One may notice the significant gap between the former lower bounds for weak Schur numbers obtained by conducting a computer search and the new lower bounds obtained with WS-templates; thus questioning the methods used and the assumptions made when obtaining these former lower bounds. We reckon better WS-templates with $n \geqslant 5$ colors can be found by making different assumptions and using a different method (Monte-Carlo methods for instance).

# 4 Conclusion and future work

We have produced new templates for Schur and Ramsey numbers. We also have generalized the concept of template to weak Schur numbers. These templates provide new general inequalities of the form $S(n+k) \geqslant aS(n)+b$ and $WS(n+k) \geqslant aS(n)+b$ as well as new lower bounds for both Schur and weak Schur numbers. We have given a new lower bound $WS(6) \geqslant 646$ with a method that can yield slight improvements for fixed values of $n$ over lower bounds obtained with the inequality $WS(n+1) \geqslant 4S(n)+2$.

Given that the introduction of templates is quite recent, we expect the bounds to be improved as this special type of partition becomes better understood and larger templates are found. One may try to produce better templates by using heuristics or approximation algorithms, such as Monte-Carlo algorithms for instance.

The best lower bound $WS(6) \geqslant 583$ achieved with a computer search using Monte-Carlo methods [4] is significantly lower than those obtained with a template-like structure: $WS(6) \geqslant 642$ [8], $WS(6) \geqslant 646$. We have evidence that 583 is the maximal value in the search space considered by [4, 5, 6]. Therefore, it questions the assumption that large partitions for $WS(n+1)$ start with large partitions for $WS(n)$. It also indicates that $WS(5)$ might need further investigation since the current lower bound $WS(5) \geqslant 196$ [12] was obtained by considering the same type of search space.

# 5 Acknowledgements

This study results from a long-term project performed within the framework of the Project Cluster "Training for Research" of CentraleSupélec. We would like to address warm thanks to Fred Rowley who kindly answered our many questions and with who it was a pleasure to exchange. Last but not least, we would like to thank Marijn Heule for providing the backdoors used in the computation of the new lower bound for $WS(6)$.

# A  S-templates

Table 5: S-template with width 33 and 4 colors

| | |
|---|---|
| 1 | $1, 6, 9, 13, 16, 20, 24, 27, 31$ |
| 2 | $2, 5, 14, 15, 25, 26$ |
| 3 | $3, 4, 10, 11, 12, 28, 29, 30$ |
| 4 | $7, 8, 17, 18, 19, 21, 22, 23, 32, 33$ |

Table 6: S-template with width 111 and 5 colors

| | |
|---|---|
| 1 | $1, 5, 8, 12, 14, 21, 23, 30, 32, 36, 39, 43, 45, 52, 103, 106, 110$ |
| 2 | $2, 6, 7, 10, 15, 18, 26, 29, 34, 37, 38, 42, 46, 51, 54, 101, 104, 109$ |
| 3 | $3, 4, 9, 11, 17, 19, 25, 27, 33, 35, 40, 41, 47, 48, 55, 100, 107, 108$ |
| 4 | $13, 16, 20, 22, 24, 28, 31, 58, 61, 67, 88, 94, 97$ |
| 5 | $44, 50, 53, 56, 57, 59, 60, 62, 63, 64, 65, 66, 68, 69, 70,$ $71, 72, 73, 74, 75, 76, 77, 78, 79, 80, 81, 82, 83, 84, 85,$ $86, 87, 89, 90, 91, 92, 93, 95, 96, 98, 99, 102, 105, 111$ |

Table 7: S-template with width 380 and 6 colors

| 1 | 1, 5, 8, 11, 15, 17, 29, 33, 36, 39, 43, 57, 61, 88, 92, 106, 110, 113, 116, 120, 132, 134, 138, 141, 144, 148, 150, 154, 157, 160, 164, 178, 182, 185, 188, 341, 344, 347, 351, 365, 369, 372, 375, 379 |
|---|---|
| 2 | 2, 9, 13, 16, 20, 23, 24, 27, 28, 31, 34, 35, 38, 42, 45, 49, 53, 60, 67, 71, 78, 82, 89, 96, 100, 104, 107, 111, 114, 115, 118, 121, 122, 125, 126, 129, 133, 136, 140, 147, 158, 162, 165, 169, 172, 176, 183, 187, 194, 201, 328, 335, 342, 346, 353, 357, 360, 364, 367, 371 |
| 3 | 3, 4, 12, 14, 19, 25, 30, 32, 40, 41, 47, 48, 58, 91, 101, 102, 108, 109, 117, 119, 124, 130, 135, 137, 145, 146, 152, 153, 161, 163, 168, 179, 181, 190, 339, 348, 350, 361, 366, 368, 376, 377 |
| 4 | 6, 7, 10, 18, 21, 22, 26, 37, 46, 50, 51, 54, 65, 70, 79, 84, 95, 98, 99, 103, 112, 123, 127, 128, 131, 139, 142, 143, 151, 155, 156, 159, 167, 170, 171, 175, 186, 343, 354, 358, 359, 362, 370, 373, 374, 378 |
| 5 | 44, 52, 55, 56, 59, 62, 63, 64, 66, 68, 69, 72, 73, 74, 75, 76, 77, 80, 81, 83, 85, 86, 87, 90, 93, 94, 97, 105, 189, 196, 197, 200, 203, 206, 207, 209, 214, 219, 231, 298, 310, 315, 320, 322, 323, 326, 329, 332, 333, 340 |
| 6 | 149, 166, 173, 174, 177, 180, 184, 191, 192, 193, 195, 198, 199, 202, 204, 205, 208, 210, 211, 212, 213, 215, 216, 217, 218, 220, 221, 222, 223, 224, 225, 226, 227, 228, 229, 230, 232, 233, 234, 235, 236, 237, 238, 239, 240, 241, 242, 243, 244, 245, 246, 247, 248, 249, 250, 251, 252, 253, 254, 255, 256, 257, 258, 259, 260, 261, 262, 263, 264, 265, 266, 267, 268, 269, 270, 271, 272, 273, 274, 275, 276, 277, 278, 279, 280, 281, 282, 283, 284, 285, 286, 287, 288, 289, 290, 291, 292, 293, 294, 295, 296, 297, 299, 300, 301, 302, 303, 304, 305, 306, 307, 308, 309, 311, 312, 313, 314, 316, 317, 318, 319, 321, 324, 325, 327, 330, 331, 334, 336, 337, 338, 345, 349, 352, 355, 356, 363, 380 |

# B  WS-templates

Table 8: 23-WS-template with width 42 and 4 colors

| 1 | 1, 2, 4, 8, 11, 22, 25, 48, 53, (**N** + **1**) |
|---|---|
| 2 | 3, 5, 6, 7, 19, 21, 23, 36, 50, 51, 52, 63, 64, 65 |
| 3 | 9, 10, 12, 13, 14, 15, 16, 17, 18, 20, 54, 55, 56, 57, 58, 59, 60, 61, 62 |
| 4 | 24, 26, 27, 28, 29, 30, 31, 32, 33, 34, 35, 37, 38, 39, 40, 41, 42, 43, 44, 45, 46, 47, 49 |

This template provides the inequality $WS(n+3) \geqslant 42\,S(n)+24$ by placing one last number, here represented by (**N** + **1**), in the first subset (see Proposition 3.20).

# C  $WS(6) \geqslant 646$

Table 9: Weakly sum-free partition of $[\![1, 646]\!]$ into 6 subsets

| 1 | 1, 2, 6, 10, 14, 18, 26, 30, 34, 42, 46, 50, 54, 62, 70, 79, 82, 90, 95, 99, 111, 115, 119, 123, 131, 135, 139, 143, 151, 155, 159, 163, 171, 175, 179, 183, 187, 195, 199, 203, 207, 211, 215, 220, 224, 228, 232, 236, 239, 244, 252, 256, 260, 264, 267, 272, 275, 280, 284, 292, 296, 300, 304, 308, 312, 316, 320, 328, 340, 344, 348, 353, 360, 364, 368, 372, 381, 385, 388, 393, 397, 404, 408, 413, 417, 425, 428, 433, 441, 445, 449, 453, 457, 461, 465, 469, 473, 485, 489, 493, 497, 502, 505, 509, 513, 517, 521, 525, 529, 533, 537, 541, 546, 549, 553, 558, 562, 566, 569, 574, 578, 586, 590, 593, 598, 602, 606, 610, 614, 618, 622, 626, 630, 634, 638, 642, 646 |
|---|---|
| 2 | 3, 4, 5, 15, 16, 17, 27, 28, 29, 39, 40, 41, 47, 48, 49, 112, 113, 114, 120, 121, 122, 132, 133, 134, 156, 157, 158, 176, 177, 178, 200, 201, 202, 221, 222, 258, 259, 281, 282, 283, 301, 302, 303, 345, 346, 347, 365, 366, 367, 389, 426, 427, 446, 447, 448, 470, 471, 472, 490, 491, 492, 514, 515, 516, 526, 527, 528, 534, 535, 536, 599, 600, 601, 607, 608, 609, 619, 620, 621, 631, 632, 633, 643, 644, 645 |
| 3 | 7, 8, 9, 19, 20, 21, 22, 23, 24, 25, 35, 36, 37, 38, 87, 88, 89, 136, 137, 138, 150, 152, 153, 154, 180, 181, 182, 196, 197, 198, 208, 209, 210, 212, 213, 214, 261, 262, 263, 265, 266, 276, 277, 278, 279, 309, 321, 322, 323, 324, 325, 326, 327, 338, 339, 369, 370, 371, 382, 384, 386, 387, 434, 435, 436, 437, 438, 439, 440, 450, 451, 452, 466, 467, 468, 482, 494, 495, 496, 499, 500, 510, 511, 512, 559, 560, 561, 611, 612, 613, 623, 624, 625, 627, 628, 629, 639, 640, 641 |
| 4 | 11, 12, 13, 31, 32, 33, 51, 100, 101, 102, 103, 104, 105, 106, 107, 108, 109, 110, 124, 125, 126, 127, 128, 129, 130, 144, 145, 146, 147, 148, 149, 164, 165, 166, 167, 168, 169, 170, 172, 173, 174, 188, 189, 190, 191, 192, 193, 194, 454, 455, 456, 458, 459, 460, 474, 475, 476, 477, 478, 479, 480, 481, 483, 484, 498, 501, 503, 504, 518, 519, 520, 522, 523, 524, 538, 539, 540, 542, 543, 544, 545, 547, 548, 597, 615, 616, 617, 635, 636, 637 |
| 5 | 43, 44, 45, 52, 53, 55, 56, 57, 58, 59, 60, 61, 63, 64, 65, 66, 67, 68, 69, 71, 72, 73, 74, 75, 76, 77, 78, 80, 81, 83, 84, 85, 86, 91, 92, 93, 94, 216, 217, 218, 219, 223, 225, 226, 227, 231, 233, 234, 235, 237, 238, 240, 241, 242, 243, 245, 246, 247, 251, 253, 254, 255, 257, 383, 390, 391, 392, 394, 395, 396, 400, 401, 402, 403, 405, 406, 407, 409, 410, 411, 412, 414, 415, 416, 421, 422, 423, 424, 429, 430, 431, 432, 554, 555, 556, 557, 563, 564, 565, 567, 568, 570, 571, 572, 573, 575, 576, 577, 579, 580, 581, 582, 583, 584, 585, 587, 588, 589, 591, 592, 594, 595, 596, 603, 604, 605 |
| 6 | 96, 97, 98, 116, 117, 118, 140, 141, 142, 160, 161, 162, 184, 185, 186, 204, 205, 206, 229, 230, 248, 249, 250, 268, 269, 270, 271, 273, 274, 285, 286, 287, 288, 289, 290, 291, 293, 294, 295, 297, 298, 299, 305, 306, 307, 310, 311, 313, 314, 315, 317, 318, 319, 329, 330, 331, 332, 333, 334, 335, 336, 337, 341, 342, 343, 349, 350, 351, 352, 354, 355, 356, 357, 358, 359, 361, 362, 363, 373, 374, 375, 376, 377, 378, 379, 380, 398, 399, 418, 419, 420, 442, 443, 444, 462, 463, 464, 486, 487, 488, 506, 507, 508, 530, 531, 532, 550, 551, 552 |